\documentstyle{article}
\def\d{\partial}
\def\dfrac#1#2{{\displaystyle {#1 \over #2}}}

\def\dsum{\mathop{\displaystyle \sum }}
\def\stackunder#1#2{\mathrel{\mathop{#2}\limits_{#1}}}
\font\bb=msbm10 at 11 pt

\def\R{\hbox{\bb R}}
\def\QED {\hfill \vrule height 1.2ex width 1.2ex}{\vskip 10pt plus100pt}

\begin{document}
\noindent Analyse Fonctionnelle/Functional Analysis
\bigskip \bigskip

\centerline{\large \bf Limite classique de C$^*$-alg\`ebres de 
groupo{\"\i}des de Lie}
\medskip

\centerline{Birant RAMAZAN}
\medskip

{\it Adresse: } Institute of Mathematics, Romanian Academy, 
\\ \hspace*{2.2cm}P.O.Box 1-764, 70700 Bucharest, Romania
\\ \hspace*{2.1cm} e-mail : ramazan@stoilow.imar.ro  
\bigskip

{\bf R\'esum\'e } \parbox[t]{12cm}{\small Soit $G$ un groupo{\"\i}de de Lie et
$\cal G$ son alg\'ebro{\"\i}de de Lie. On propose une d\'efinition de limite 
classique d'un champ continu de C$^*$-alg\`ebres et on prouve que le 
groupo{\"\i}de tangent $\tilde{G}$ associ\'e \`a $G$ permet d'obtenir la 
structure canonique de Poisson sur $\cal G$ 
comme limite classique d'un champ continu de C$^*$-alg\`ebres.}
\medskip

\centerline{\bf Classical limit of Lie groupoid C$^*$-algebras}
\smallskip

{\bf Abstract } \parbox[t]{12cm}{\small Let $G$ be a Lie groupoid and 
${\cal G}$ his Lie algebroid. We give a definition of the classical limit 
of a C$^*$-bundle and we use the tangent groupoid $\tilde{G}$ associated to 
$G$ to show that the Poisson structure on ${\cal G}$ is the classical limit 
of a C$^*$-bundle.} 
\bigskip

{\it 1991 Mathematics Subject Classification.} Primary 46L60; Secondary 46L87,
81S99, 22A22

\bigskip \bigskip

On propose une d\'efinition de limite classique en termes de C$^*$-alg\`ebres.

\medskip

\noindent {\bf D\'efinition 1} {\em Soit $(A^0,\left\{ \cdot ,\cdot \right\} )
$ une $*$-alg\`ebre de Poisson, $I$ un intervalle dans $\R$ contenant $0$, et 
$(A_t)_{t\in I}$ un champ continu de C$^*$-alg\`ebres d'espace total $A$. 
On dira que $(A^0,\left\{ \cdot ,\cdot \right\} )$
est la limite classique du champ $(A_t)_{t\in I}$ si :
\begin{enumerate}
\item $A^0$ est une sous-$*$-alg\`ebre dense de $A_0$
\item Il existe une sous-alg\`ebre dense de sections ${\cal A}\subset C_0(I 
, A)$ telle que :
\begin{enumerate}
\item $\left\{ f(0), f\in {\cal A}\right\}$ est un sous-ensemble dense de $A^0$
\item l'application $z:{\cal A}\times {\cal A}\rightarrow {\cal A}$, donn\'ee
par $z(f,g)(t)=\dfrac{ \left[ f(t),g(t)\right] }{it}$, pour $t\neq 0$ et 
$z(f,g)(0)=\left\{ f(0),g(0)\right\}$, d\'efinit une structure d'alg\`ebre 
de Poisson sur ${\cal A}$.
\end{enumerate}
\end{enumerate}}

\bigskip

{\it Notations}

Soit $G$ un groupo{\"\i}de de Lie de dimension $n+m$, $G^{(0)}$ son espace 
d'unit\'es de dimension $n$, et
$r$, respectivement $s$, l'application but, respectivement source de $G$. 
On consid\`ere l'alg\`ebre involutive de convolution $C_c^{\infty}(G, 
\Omega ^{1/2}(G))$ des demi-densit\'es sur $G$, alg\`ebre d\'efinie par 
A.~Connes dans \cite{connes}.
On note C$^*(G)$ la C$^*$-alg\`ebre enveloppante de $C_c^{\infty}(G,
\Omega ^{1/2}(G))$ et C$^*_r(G)$ la C$^*$-alg\`ebre r\'eduite. 

Soit aussi ${\cal G}$ l'alg\'ebro{\"\i}de de Lie associ\'e \`a $G$,
introduit par J.~Pradines dans \cite{pradines}. Il existe
une bijection entre les sections de $\Omega ^{1/2}({\cal G})$ et
les sections \'equivariantes \`a gauche de $\Omega ^{1/2}(G)$. 
On fixe $\mu$ une section $C^{\infty}$ de $\Omega ^{1/2}({\cal G})$ et on lui 
associe la section \'equivariante \`a gauche $\lambda$ de $\Omega ^{1/2}(G)$
(ces deux sections joueront le r\^ole d'un syst\`eme de Haar pour ${\cal G}$,
respectivement $G$). Cela  
permet de se r\'eduire aux fonctions dans les calculs, et d'utiliser
les r\'esultats de la th\'eorie des C$^*$-alg\`ebres de groupo{\"\i}des 
d\'evelopp\'ee par J.~Renault (\cite{renault}). 

\medskip

Le r\'esultat qu'on prouvera dans cette note, li\'e \`a l'id\'ee 
de Landsman (cf. \cite{landsman}) de la d\'eformation d'un 
alg\'ebro{\"\i}de de Lie par le groupo{\"\i}de de Lie correspondant, 
est le suivant :

\noindent {\bf Th\'eor\`eme 1} 
{\em La structure de Poisson canonique sur l'alg\`ebre de convolution 
$S({\cal G})$ des fonctions de type Schwartz, est la 
limite classique d'un champ continu sur $\R$ de C$^*$-alg\`ebres, ayant 
C$^*({\cal G})$ comme fibre en $t=0$ et C$^*(G)$ au dessus des $t\neq 0$.}
\bigskip

\centerline{\bf 1. L'alg\`ebre de Poisson $S({\cal G})$}
\medskip

Le dual ${\cal G}^*$ de ${\cal G}$ poss\`ede une structure canonique de 
vari\'et\'e de Poisson, mise en \'evidence par A.~Weinstein, P.~Dazord et 
A.~Coste \cite{coste}. Cette structure peut \^etre transport\'ee sur ${\cal G}$
en utilisant la transformation de Fourier qui est un isomorphisme de 
C$^*$-alg\`ebres ${\cal F}:C^*({\cal G})\rightarrow C_0({\cal G}^*)$.

On introduit l'espace de Schwartz $S(E)$ d'un fibr\'e vectoriel $E$, en 
g\'en\'eralisant la d\'efinition introduite par M.~Rieffel
dans le cas du fibr\'e trivial $M\times \R ^d$. La restriction de la 
transformation de Fourier ${\cal F}$ est un isomorphisme d'alg\`ebres entre 
$(S({\cal G}),+,*)$ et $(S({\cal G}^*),+,\cdot )$ et on obtient un crochet de 
Poisson sur l'alg\`ebre de convolution $S({\cal G})$ par 
$\left\{ f,g\right\} ={\cal F}^{-1}(
\left\{ {\cal F}f,{\cal F}g\right\} )$.

Par calcul direct on prouve que si $f,g\in S({\cal G})$ alors le crochet de 
Poisson de $f$ et $g$ est donn\'e localement par  
\begin{eqnarray}\label{crosetE}
\left\{ f,g\right\} (x,\xi ) & = & 2\pi i\stackunder{i,j}{\dsum}
a_{ij}(x) \left( \xi _if*\dfrac{\d g}{\d x_j} -\xi _ig*\dfrac{\d f}{\d x_j}
\right) (x,\xi ) \nonumber \\
 & +& 2\pi i \stackunder{i,j}{\dsum}a_{ij}(x)\dfrac{\d \ln \mu _e}{\d x_j}(x)
(\xi _if*g -\xi _ig*f)(x,\xi ) \nonumber \\
 & -& 2\pi i\stackunder{i,j,k}{\dsum}c_{ijk}(x)\dfrac{\d }{\d \xi _k}
\left( \xi _if*\xi _jg\right) (x,\xi )
\end{eqnarray}

\noindent o\`u $(e_1,..,e_m)$ est un rep\`ere mobile de ${\cal G}$ et 
$\mu _e$ est la fonction donn\'ee par $\mu _e(x)=\mu (x)(e_1(x)\wedge ..\wedge
e_m(x))$. 
\bigskip

\centerline{\bf 2. La structure locale de $G$}

\medskip

On explicite localement la structure de $G$ dans une carte choisie 
convenablement au voisinage de $x_0\in G^{(0)}\subset G$, et on montre comment 
on lui associe une carte de ${\cal G}$.
Les cartes de cette forme ont \'et\'e consid\'er\'ees aussi par 
Weinstein, Nistor et Xu dans \cite{nistor}.

Soit $x_0$ appartenant \`a la sous-vari\'et\'e $G^{(0)}$ de $G$. 
Comme $r$ est une submersion dans le point $x_0$, il existe $U$ 
voisinage ouvert de 0 dans 
$\R^n$, $V$ voisinage ouvert de 0 dans $\R^m$ et les cartes $\psi :U\times 
V\longrightarrow G$, $\varphi :U\longrightarrow G^{(0)}$ v\'erifiant:

\begin{equation}\label{psibon}
\begin{array}{l}
1.\ \psi (0,0)=x_0\\
2.\ r(\psi (u,v))=\varphi (u)\\
3.\ \psi (U\times \lbrace 0\rbrace )=\psi (U\times V)\bigcap G^{(0)}
\end{array}
\end{equation}

\noindent A la carte $\psi$ de $G$ s'associe la carte $\theta :U\times \R^m
\longrightarrow {\cal G}$, $\theta (u,v)=(\varphi (u),\dfrac{\d \psi}{\d v}
(u,0)v)$ de ${\cal G}$ au voisinage de la fibre ${\cal G}_{x_0}$. La famille 
$\lbrace e_1,e_2,..,e_m\rbrace$, $e_i(\varphi (u))=\theta (u,f_i)$, 
$i=\overline{1,\ldots,m}$,
o\`u $\lbrace f_1,f_2,...,f_m\rbrace $ est la base canonique de $\R^m$, est
un rep\`ere mobile de ${\cal G}$ sur $\varphi (U)$.

Notons $\sigma :U\times V\longrightarrow U$, $s(\psi (u,v))=\varphi (\sigma 
(u,v))$ la submersion qui d\'ecrit l'application source $s$ dans 
les cartes $\psi$ et $\varphi$. On a $\sigma (u,0)=u$.
  
Le r\'esultat suivant donne une formule de type Baker-Campbell-Hausdorff 
pour $G$.
\newpage
\noindent {\bf Proposition 2} 

{\em (i) $(\psi (u,v),\psi (u_1, w))\in G^{(2)}$ si et seulement si 
$u_1=\sigma (u,v)$. Dans ce cas leur produit est donn\'e par
\begin{equation}\label{produs}
\psi (u,v)\psi (\sigma (u,v),w)=\psi (u,p(u,v,w))
\end{equation}

\noindent o\`u $p:U\times V\times V\rightarrow V$ est une application
diff\'erentiable qui a un d\'eveloppement de la forme 

\centerline{$p(u,v,w)=v+w+B(u,v,w)+O_3(u,v,w)$}

\smallskip

\noindent avec $B$ bilin\'eaire en $(v,w)$ et $O_3(u,v,w)$ 
de l'ordre de $\|(v,w)\|^3$.

(ii) Soit $(u,v)\in U\times V$ tel que $\psi (u,v)^{-1}\in \psi (U\times 
V)$. Alors $\psi (u,v)^{-1}=\psi (\sigma (u,v),w)$, o\`u $w$ v\'erifie 
$p(u,v,w)=0$. De plus on a le d\'eveloppement $w=-v+B(u,v,v)+O_3(u,v)$, avec
$O_3(u,v)$ de degr\'e d'homog\'en\'eit\'e superieur \`a 3 en $v$.}

\medskip

\noindent {\bf Preuve} Soit $g=\psi (u,v)$ et $h=\psi (u_1,v)$. On a 
$s(g)=\varphi (\sigma (u,v))$ et $r(h)=\varphi (u_1)$ ce qui montre la 
premi\`ere assertion. De plus $r(gh)=r(g)=\varphi (u)$ assure l'existence d'un
unique $p(u,v,w)\in V$ tel qu'on a l'\'egalit\'e \ref{produs}. L'application
$p$ est bien d\'efinie et v\'erifie $p(u,0,w)=w$, $p(u,v,0)=v$. On d\'erive 
ces deux \'egalit\'es et on les utilise dans le d\'eveloppement de Taylor de 
$p$.

\QED

\medskip

\noindent {\bf Remarque} On peut calculer les fonctions de structure de l'alg\'ebro{\"\i}de 
${\cal G}$ qui pour tout $u\in U$
sont donn\'ees par $a_{ij}(\varphi (u))=\dfrac{\d \sigma _j}{\d 
v_i}(u,0)$ et $c_{ijk}(\varphi (u))=B_k(u,f_i,f_j)-B_k(u,f_j,f_i)$.

\bigskip

\centerline{\bf 3. Le groupo{\"\i}de tangent $\tilde{G}$}

\medskip

Le groupo{\"\i}de tangent de $G$ introduit par M.~Hilsum et G.~Skandalis 
(\cite{hilsum}), est le groupo{\"\i}de de Lie $\tilde{G}=G\times \R ^*
\bigcup {\cal G}\times \left\{ 0\right\}$ qui a la structure alg\'ebrique 
de fibr\'e des groupo{\"\i}des $G(t)=G\times \left\{ t\right\}$, pour $t\neq 0$
et $G(0)={\cal G}\times \left\{ 0\right\}$ et la structure de vari\'et\'e
obtenue comme \'eclatement de $G$ au long de $G^{(0)}$.
Pour $G=M\times M$ on retrouve le groupo{\"\i}de tangent de A.~Connes 
(\cite{connes}).

\bigskip

\noindent {\bf Proposition 3} {\em Les C$^*$-alg\`ebres C$^*(G(t))$, 
$t\in \R$, forment un champ continu dont la C$^*$-alg\`ebre est  
C$^*(\tilde{G})$.}

\medskip

\noindent {\bf Preuve :} Les fonctions de $C_c(\tilde G)$ d\'efinissent par 
restrictions une famille des sections du fibr\'e des C$^*(G(t))$ sur $\R$. 
Cette famille est stable pour la convolution et l'involution, et de plus 
pour chaque $t\in \R$, $C_c(G(t))$ est dense dans C$^*(G(t))$.

Soit $f\in C_c(\tilde G)$. Puisque le fibr\'e est trivial dans les points
$t\neq 0$, l'application $t\mapsto \| f_t\|$ est continue dans ces points. 
La continuit\'e en $t=0$ est assur\'ee par le th\'eor\`eme suivant, car
${\cal G}$ est commutatif, donc  moyennable. 

\QED
\bigskip

\noindent {\bf Th\'eor\`eme 4}
{\em Soit $G\stackrel{p}{\longrightarrow}X$ un champ continu de groupo{\"\i}des et
$a\in C_c(G)$. Pour $x\in X$, on note $a_x=a_{|_{G(x)}}$. Alors :

i) L'application $X\ni x\mapsto \left\| a_x\right\| _{C^*(G(x))}$ est 
semi-continue sup\'erieurement.

ii) L'application $X\ni x\mapsto \left\| a_x\right\| _{C^*_r(G(x))}$ est
semi-continue inf\'erieurement.}
\bigskip

\noindent {\bf Preuve } La d\'emonstration est bas\'e sur le th\'eor\`eme de 
d\'esint\'egration des repr\'esentations de J.~Renault (\cite{renault}) et les 
travaux de E.~Blanchard (\cite{blanchard}) sur les $C(X)$-alg\`ebres 
stellaires.

\QED

\bigskip

Le calcul principal, n\'ecessaire dans la d\'emonstration du th\'eor\`eme 1 
est contenu dans le lemme technique suivant : 

\medskip

\noindent {\bf Lemme 5}
{\em Soit $f,g\in C_c^{\infty }(\tilde{G})$ et $f_0=f|_{\tilde{G}(0)}$, 
$g_0=g|_{\tilde{G}(0)}$. Alors 
$$\dfrac{\d }{\d t}\lbrack f,g\rbrack _{|_{t=0}}=\dfrac{1}{2\pi i}
\lbrace f_0,g_0\rbrace$$}
\bigskip

\noindent {\bf Preuve } On d\'efinit de mani\`ere \'evidente un champ de 
vecteurs $\dfrac{\d}{\d t}\in C^{\infty}(\tilde{G},T\tilde{G})$.

Par $\tilde{\lambda}(g,t)=\dfrac{1}{|t|^m}\lambda (g)$,
pour $g\in G$, et $\tilde{\lambda} (x,X,0)=\lambda (x)$, pour 
$(x,X)\in {\cal G}$ 
on associe \`a $\lambda$ une section \'equivariante \`a gauche de $\Omega 
^{1/2}(\tilde{G})$, section qu'on utilise comme syst\`eme de Haar dans le 
calcul de la convolution de $f$ et $g$. On utilise la proposition 2,
on exprime le fait que $\lambda$ est \'equivariante \`a gauche et par 
calcul effectif on retrouve le crochet de Poisson donn\'e par \
la formule (\ref{crosetE}).

\QED

\bigskip

\noindent {\bf Preuve du th\'eor\`eme 1} : 
La continuit\'e du champ est prouv\'ee dans la proposition 3. 

On consid\`ere ${\cal A}=C_c^{\infty}(\tilde{G})$. On a 
$\left\{ f(0), f\in {\cal A}\right\} =C_c^{\infty}({\cal G})$. 
Les conditions 1. et 2.(a) de la d\'efinition 1 sont cons\'equences 
du fait que $S({\cal G})$ contient la sous-alg\`ebre 
$C_c^{\infty}({\cal G})$ dense dans $C^*({\cal G})$.

Enfin la condition 2.(b) est la cons\'equence imm\'ediate du lemme 5. 

\QED

\bigskip

\noindent {\it Remerciemements.} L'auteur tient a exprimer sa gratitude 
\`a Jean Renault et Jean Pradines pour leurs judicieuses remarques et au 
D\'epartement de Math\'ematiques de l'Universit\'e d'Orl\'eans o\`u ce 
travail a \'et\'e effectu\'e.

\end{document}